\documentclass[a4paper,12pt]{article}
\usepackage{amsmath}
\usepackage{amsfonts}
\usepackage{amssymb}
\usepackage{latexsym}
\usepackage{epsfig}
\usepackage{graphicx}
\usepackage{oldgerm}
\usepackage{theorem}
\usepackage{graphicx}

\makeatletter
 
 \@addtoreset{equation}{section}
\makeatother

\def\N{\mathbb{N}}

\def\R{\mathbb{R}}

\def\Z{\mathbb{Z}}

\def\S{\widetilde{S}}
\def\SL{\widetilde{S}^{f-}}
\def\SR{\widetilde{S}^{f+}}
\def\Lf{L^{f}}
\def\Rf{R^{f}}
\def\h{h_{\delta}}

\def\0{\mib{0}}

\def\={\stackrel{\rm (law)}{=}}

\theorembodyfont{\itshape}

\newtheorem{thm}{Theorem}
\newtheorem{lem}[thm]{Lemma}
\newtheorem{cor}[thm]{Corollary}
\newtheorem{prop}[thm]{Proposition}

\newtheorem{remark}{Remark}

\newcommand{\mib}[1]{\mbox{\boldmath $#1$}}
\newcommand{\qed}{\hbox{\rule[-2pt]{3pt}{6pt}}}

\title{On a model of evolution of subspecies}
\author{{
Rahul Roy and Hideki Tanemura}
\footnote{E-Mail: {\tt rahul@isid.ac.in} and {\tt tanemura@math.keio.ac.jp}}\\
{\it Indian Statistical Institute, New Delhi and Keio University, Yokohama}}

\date{}

\begin{document}

\maketitle

\begin{abstract}
Ben-Ari and Schinazi (2016) introduced a stochastic model to study `virus-like evolving population with high mutation rate'. This model is a birth and death model with an individual at birth being either a mutant with a random fitness parameter in $[0,1]$ or having one of the existing fitness parameters with uniform probability; whereas a death event removes the entire population of the least fit site. We change this to incorporate the notion of `survival of the fittest', by requiring  that a non-mutant individual, at birth, has a fitness according to a preferential attachment mechanism, i.e., it has a fitness $f$ with a probability proportional to the size of the population of fitness $f$. Also death just removes one individual at the least fit site. This preferential attachment rule leads to a power law behaviour in the asymptotics, unlike the exponential behaviour obtained by Ben-Ari and Schinazi (2016).
\end{abstract}

\vspace{0.1in}

\vspace{0.1in}

\section{Introduction}
We study a model of the evolution and survival of species subjected to birth, mutation and death. This model was introduced by Guiol, Machado and Schinazi (2010) and is similar to a model studied by Liggett and Schinazi (2009). It has been of recent interest because of its relation to the discrete evolution  model of Bak and Sneppen (1993). 

In the model studied by Guiol, Machado and Schinazi (2010), at each discrete time point,  with probability $p$ or $1-p$ respectively,  there is either a birth of an individual of the species or a death (in case there exists at least one surviving species). An individual at birth is accompanied by a fitness parameter $f$, which is chosen uniformly in $[0,1]$, while the death is always of the individual with the least fitness parameter. They exhibited a phase transition in this model, i.e., for $p > 1/2$, the size of the population, $L_n$, at time $n$ whose fitness is smaller that $f_c := (1-p)/p$ is a null recurrent Markov chain, while asymptotically, the proportion of the population with fitness level  lying in $(a,b)\subseteq (f_c, 1)$ equals $p(b-a)$ almost surely.

In a subsequent paper Ben-Ari and Schinazi (2016) modified the above model to study a `virus-like evolving population with high mutation rate'. Here, as earlier, at each discrete time point, with probability $p$ or $1-p$ respectively, there is either a birth of an individual of the species or a death (in case there exists at least one surviving species) of the individual with the least fitness parameter. 
The caveat here is that at death, the entire population of the least fit individuals is removed; while, at birth, the individual,
\begin{itemize}
\item[(i)]  with probability $r$, is a mutant and has a fitness parameter $f$ uniformly at random in $[0,1]$,  or 
\item[(ii)] with  probability $1-r$, has a fitness parameter chosen uniformly at random among the existing fitness parameters,  thereby increasing the population at that fitness level by $1$. 
\end{itemize}
For this model too, the authors exhibited a phase transition. In particular, assuming $pr > (1-p)$, for $f_c := (1-p)/pr$ the number of fitness levels lying in $(0, f_c)$ at time $n$ where individuals exist  is a null recurrent Markov chain, while the number of fitness levels lying to the right of $f_c$  is asymptotically uniformly distributed in $(f_c, 1)$ uniformly. 

Here we propose a variant of the Ben-Ari, Schinazi model, a variant which we believe is closer to the Darwinian theory of the survival of the fittest. 
To incorporate the Darwinian theory, we differ from the above model when  a birth occurs which is not a mutant. Instead of the individual at birth having a fitness one of the  existing fitness levels chosen uniformly at random, the newly born individual has a fitness $f$ which is chosen proportional to the size of the population of fitness $f$.

 More particularly, suppose that at time $n$ there is a birth, which is not a mutant, and that there are $n_i$ individuals with fitness $f_i$ for $i = 1, \ldots , k$ and no other individuals elsewhere. The newly born individual has a fitness $f_j$ with a probability proportional to ${n_j}$ for $j = 1, \ldots, k$. Thus, at birth, an individual without mutation follows a preferential attachment rule akin to the  Barab\'{a}si and Albert (1999) model.
 
Before we end this section we note that Schreiber (2001) and subsequently Bena\"{i}m, 
Schreiber and Tarr\`{e}s (2004) study the question of random genetic drift and natural selection via urn models coupled with mean-field behaviour. Unlike our study, there is no spatial aspect of fitness in their model.

A formal set-up of this model is given in the next section, while in the last section we present some mean-field dynamics of the model.

\section{The model and statement of results}
We first present our model and state the results.

At time $0$ there is one individual at site $0$. At time $n$, there is either a birth or a death of an individual from the existing population with probability $p$ or $1-p$ respectively, where $p\in (0,1)$, and independent of any other random mechanism considered earlier.
\begin{enumerate}
\item[(P1)] In case of a birth, there are two possibilities. 
\begin{itemize}
\item[(i)]  with probability $r\in (0,1)$,  a mutant is born and has a fitness parameter $f$ uniformly at random in $[0,1]$,  or 
\item[(ii)] with 
probability $1-r$ the individual born has a fitness $f$ with a probability proportional to the number of individuals with fitness $f$ among the entire population present at that time. Here we have a caveat that, if there is no individual present at the time of birth, then the fitness of the individual is  sampled uniformly in  $[0,1]$.
\end{itemize}
\item[(P2)] In case of a death, an individual from the population at the  site closest to $0$ is eliminated.
\end{enumerate}
Here and henceforth, a site represents a fitness level.

Let $X_n = \{(k_i, x_i) : k_i \geq 1, x_i \in [0,1], i = 1, \ldots , l\}$, where the total population at time $n$ is divided in exactly $l$ sites $x_1, \ldots , x_l$, with the size of the population at site $x_i$ being exactly $k_i$. In case there is no individual present at time $n$ we take $X_n = \emptyset$.
The process $X_n$ is Markovian on the state space
\begin{align}
\mathbb{X}:=\{\emptyset\} \cup \{ \{ (k,x)\}_{x\in \Lambda} : (k,x)\in \N \times [0,1], \; \sharp \Lambda <\infty, \},
\end{align}
where $\N = \{1,2,\dots\}$.

For a given $f\in (0,1)$, let $L_n^f$ denote the size of the population at time $n$ at sites in [0,f],
$$
L_n^f :=  \sum_{s \in [0,f]}k_s : s\in [0,f] \text{ and }  (k_s ,s) \in X_n , 
$$
$R_n^f$  denote the size of the population at time $n$ at sites in $(f,1]$,
$$
R_n^f :=  \sum_{s \in (f,1]}k_s : s\in (f,1] \text{ and }  (k_s ,s) \in X_n, 
$$
and $N_n$  denote the size of the population at time $n$, 
$$
N_n:= L_n^f+R_n^f.
$$

For a fixed $f \in (0,1)$, the pair $(L_n^f, R_n^f)$ is a Markov chain on $\Z_+\times\Z_+$, ($\Z_+=\N \cup \{0\}$) with transition probabilities given by 

\vskip 3mm

(1-1) If $(L_n^f, R_n^f)=(0,0)$
\begin{equation}\label{TP11}
(L_{n+1}^f, R_{n+1}^f)=
\begin{cases}
(1,0) \quad &\mbox{w. p. $fp$}
\\
(0,1) &\mbox{w. p. $(1-f)p$}
\\
(0,0) &\mbox{w. p. $1-p$}
\end{cases}
\end{equation}

(1-2) If $(L_n^f, R_n^f)\in \{0\}\times \N$
\begin{equation}\label{TP2}
(L_{n+1}^f, R_{n+1}^f)=
\begin{cases}
(1,R_n^f) \quad &\mbox{w. p. $fpr$}
\\
(0,R_n^f+1) &\mbox{w. p. $(1-f)pr+p(1-r)$}
\\
(0,R_n^f-1) &\mbox{w. p. $1-p$}
\end{cases}
\end{equation}

(1-3) If $(L_n^f, R_n^f)\in \N\times \{0\}$
\begin{equation}\label{TP3}
(L_{n+1}^f, R_{n+1}^f)=
\begin{cases}
(L_n^f+1,0) \quad &\mbox{w. p. $fpr+p(1-r)$}
\\
(L_n^f,1) &\mbox{w. p. $(1-f)pr$}
\\
(L_n^f-1,0) &\mbox{w. p. $1-p$}
\end{cases}
\end{equation}

(1-4) If $(L_n^f, R_n^f)\in \N\times \N$
\begin{equation}\label{TP4}
(L_{n+1}^f, R_{n+1}^f)=
\begin{cases}
(L_n^f+1,R_n^f) \quad &\mbox{w. p. $\displaystyle{fpr+p(1-r)\frac{L_n^f}{N_n}}$}
\\
(L_n^f,R_n^f+1) &\mbox{w. p. $\displaystyle{(1-f)pr+p(1-r)\frac{R_n^f}{N_n}}$}
\\
(L_n^f-1,R_n^f) &\mbox{w. p. $1-p$}.
\end{cases}
\end{equation}
The model exhibits a phase transition at a critical position $f_c$ defined as
\begin{equation}\label{fc}
f_c:=\frac{1-p}{pr}
\end{equation}
as given in the following theorem:
\begin{thm}\label{Th1}
\begin{enumerate}
\item[\rm (1)] In case $p\le 1-p$, the population dies out infinitely often a.s., in the sense that
\begin{align}
P(N_n = 0 \text{ for infinitely many } n) = 1
\end{align}
\item[\rm \rm (2)] In case $1-p < rp$,  the size of the population goes to infinity  as $n\to\infty$, and most of the population is distributed at sites  in the interval $[f_c,1]$, in the sense that
\begin{align}
\label{Thm_12}
P(\lim_{n \to \infty} \frac{R_n^{f_c}}{N_n}  = 1) = 1 \text{ and }
P(\liminf_{n \to \infty} \frac{R_n^{f_c} - R_n^f}{N_n}  > 0) &= 1 \text{ for any } f > f_c.
\end{align}
\item[\rm (3)] In case $rp \le 1-p < p$, the size of the population goes to infinity as $n\to\infty$, and most of the population is concentrated  at sites near $1$, in the sense that
\begin{align}
\label{Thm11}
P(\lim_{n\to \infty} N_n = \infty) = 1 \text{ and, for any }\varepsilon > 0, \; P(\lim_{n \to \infty} \frac{R_n^{1-\varepsilon}}{N_n} = 1) = 1.
\end{align}

\end{enumerate}
\end{thm}

Let $F_n(f)$ denote the empirical distribution of sites at time $n$, i.e.
$$
F_n(f) := \frac{\sharp \{s \in [0,f]: (k,s) \in X_n \text{ for some } k \geq 1\}}{\sharp\{s \in[0,1]: (k,s) \in X_n \text{ for some } k \geq 1\}},
$$ 
 we have
\begin{cor}
\label{GC}
If $1-p < rp$ (i.e., $f_c < 1$), then
\begin{equation}
\label{GL1}
F_n(f) \to \frac{\max \{f-f_c, 0\}}{1-f_c}
\quad\mbox{ uniformly a.s.}
\end{equation}
\end{cor}

Let $S_n:= \sharp\{s \in[0,1]: (k,s) \in X_n \text{ for some } k \geq 1\}$ be the total number of sites at time $n$ among which  the total population is distributed.
For a given $n,k,f$ let 
$U_n^k(f):= \sharp\{s \in [f,1]: (k,s) \in X_n\}$ 
denote the number of sites in $[f,1]$ at time $n$ which has a population of size exactly $k$; clearly $S_n=\sum_{k}U_n^k(0)$. 
Taking $U_n^k(f+) = \lim_{s \downarrow f} U_n^k(s)$, for $A \subseteq \mathbb{X}$,   define the empirical distribution of size and fitness on $\N \times [0,1]$ as
\begin{equation}
H_n(A):=
\begin{cases}
\frac{\sum_{(k,f)\in A}U_n^k(f)-U_n^k(f+)}{S_n},
&S_n>0,
\\
\delta_{(0,0)}(A),
&S_n=0.
\end{cases}
\end{equation}

\begin{thm}\label{Th2}
For $pr > 1-p$, as $n \to \infty$, $H_n$ converges weakly to a product measure on $\N\times [0,1]$ whose density is given by
\begin{align}
&p_k\frac{\mathbf{1}_{[f_c, 1]}(x)}{1-f_c}dx,
\quad (k,x)\in\N\times [0,1]
\nonumber
\\
&\text{ with }
p_k=\frac{(2p-1)r}{(1-r)(1-p)} B\left(1+\frac{(2p-1)r}{(1-r)(1-p)}, k\right) \text{ for } k\in\N,
\end{align}
where $B(a,b)$ is the Beta function with parameter $a,b >0$.
\end{thm}

\begin{remark}
\label{powerlaw}
Since $B(s,k) = \mathcal{O}(k^{-s})$, $k\to\infty$,
the probability density $p_k$, $k\in\N$ has $m$-th moment if and only if $r > 1- \frac{2p-1}{2p-1 +(1-p)m}$.
\end{remark}


%
%

For the model studied by Ben-Ari and Schinazi (2016), in case of a death, the entire population at the site of lowest fitness is removed unlike our condition (P2). Thus in their model, if $\S_n$ denotes  the number of sites at time $n$ among which the total population is distributed, then 
$\S_n$ is a Markov chain with spatially homogeneous transition probabilities given by
\begin{equation}
\label{r-bas}
\S_{n+1}=\begin{cases}
\S_n+1 \quad &\mbox{with probability $pr$,}
\\
\S_n &\mbox{with probability $p(1-r)$,}
\\
\S_n-1 &\mbox{with probability $1-p$,}
\end{cases}
\end{equation}
with reflecting boundary condition at $0$.
For a given $f\in (0,1)$, letting $\SL_n$  denote the number of sites at time $n$ in [0,f],
and  $\SR_n$ the number of sites at the sites in $(f,1]$,
the pair $(\SL_n, \SR_n)$ is a spatially homogeneous Markov chain on $\Z_+\times\Z_+$, where $\Z_+=\{0,1,2,\dots \}$:

(BAS-1) If $(\SL_n, \SR_n)=(0,0)$
\begin{equation}\label{BAS1}
(\SL_{n+1}, \SR_{n+1})=
\begin{cases}
(1,0) \quad &\mbox{w. p. $fp$}
\\
(0,1) &\mbox{w. p. $(1-f)p$}
\\
(0,0) &\mbox{w. p. $1-p$}
\end{cases}
\end{equation}
(BAS-2) If $(\SL_n, \SR_n)\in \{0\}\times \N$
\begin{equation}\label{BAS2}
(\SL_{n+1}, \SR_{n+1})=
\begin{cases}
(1,\SR_n) \quad &\mbox{w. p. $fpr$}
\\
(0,\SR_n+1) &\mbox{w. p. $(1-f)pr$}
\\
(0,\SR_n) &\mbox{w. p. $p(1-r)$}
\\
(0,\SR_n-1) &\mbox{w. p. $1-p$}
\end{cases}
\end{equation}

(BAS-3) If $(\SL_n, \SR_n)\in \N\times \{0\}$
\begin{equation}\label{BAS3}
(\SL_{n+1}, \SR_{n+1})=
\begin{cases}
(\SL_n+1,0) \quad &\mbox{w. p. $fpr$}
\\
(\SL_n,1) &\mbox{w. p. $(1-f)pr$}
\\
(\SL_n,0) \quad &\mbox{w. p. $p(1-r)$}
\\
(\SL_n-1,0) &\mbox{w. p. $1-p$}
\end{cases}
\end{equation}

(BAS-4) If $(\SL_n, \SR_n)\in \N\times\N$
\begin{equation}\label{BAS4}
(\SL_{n+1}, \SR_{n+1})=
\begin{cases}
(\SL_n+1,\SR_n) \quad &\mbox{w. p. $fpr$}
\\
(\SL_n,\SR_n) &\mbox{w. p. $p(1-r)$}
\\
(\SL_n, \SR_n+1) &\mbox{w. p. $(1-f)pr$}
\\
(\SL_n-1,\SR_n) &\mbox{w. p. $1-p$}
\end{cases}
\end{equation}

Also at birth, if the individual is not a mutant then the individual born has a fitness chosen uniformly at random among the fitnesses of the existing individuals at that time, unlike the preferential condition (P1)(ii) of our model. As such, the transition probabilities for this model are spatially homogeneous, while for our model, as is exemplified by (\ref{TP4}), the transition probabilities are not spatially homogeneous. Thus the equivalent result they have for Theorem \ref{Th2} has $p_k$ arising from a $\mathrm{Geom}\left(\frac{pr-(1-p)}{p-(1-p)}\right)$ distribution.

The power law phenomenon present in the study of preferential attachment graphs (see van der Hofstad (2017) Chapter 8) manifests itself in our model (as noted in Remark \ref{powerlaw}) through the Beta function in Theorem \ref{Th2}.

%

\section{Proof of Theorem \ref{Th1}}



As noted in  Guiol, Machado and Schinazi (2010), for $p \leq 1-p$, i.e. when the death rate is more than the birth rate, the process $\{N_n : n\geq 0\}$ is equivalent to a random walk on the non-negative integers $\Z_+$ with non-positive drift and a holding at $0$ with probability $(1-p)$. Thus $N_n$ returns to the $0$ infinitely often with probability $1$. 

For $p > 1-p$, $\{N_n : n\geq 0\}$ is equivalent to a random walk on the non-negative integers $\Z_+$ with positive drift and thus $N_n \to \infty$ as $n \to \infty$ with probability $1$. 

Then we study the case when $1-p < p$.

\begin{lem}\label{lemma4}
{\rm (1)}\quad Let $f_c=\frac{1-p}{rp}<1$.  
\\
{\rm (i)}
For $f< f_c$ and for any $\eta\in (0,1)$ we have
\begin{equation}\label{key1}
P\left(\text{there exists }T>0 \mbox{ such that } \rho_n^f\equiv \frac{L_n^f}{N_n} \le \eta  \text{ for all }n\ge T \right)=1,
\end{equation}
and 
\begin{align}\label{key2}
P(L_n^f=0 \text{ infinitely often}) =1.
\end{align} 
\noindent {\rm (ii)} Let $f> f_c$.
Then 
\begin{align}\label{key3}
P(L_n^f=0 \text{ infinitely often}) = 0.
\end{align}

\noindent {\rm (2)}\quad Let $1\le f_c=\frac{1-p}{rp}<\frac{1}{r}$.  
\\
{\rm (i)}
For $f< 1$ and for any $\eta\in (0,1)$ we have (\ref{key1}) 
and (\ref{key2}).

\noindent {\rm (ii)} Let $f=1$.
Then we have $(\ref{key3})$.
\end{lem}

\noindent {\it Proof.} 
We prove two cases (1) and (2) together.
The idea of the proof is that, since for $f < f_c\wedge 1$, $R_n^f$ will be much larger than $L_n^f$, we stochastically bound the non-spatially homogeneous Markov chain with a boundary condition  by a spatially homogeneous  Markov chain a boundary condition, and study the modified Markov chain. As such, 
for $\varepsilon\in [0,1]$, we introduce a Markov chain $(\Lf_n(\varepsilon), \Rf_n(\varepsilon))$  with stationary transition probabilities given by

\vskip 3mm

(Ep-1) If $(\Lf_n(\varepsilon), \Rf_n(\varepsilon))=(0,0)$
\begin{equation}\label{TP1}
(\Lf_{n+1}(\varepsilon), \Rf_{n+1}(\varepsilon))=
\begin{cases}
(1,0) \quad &\mbox{w. p. $fp$}
\\
(0,1) &\mbox{w. p. $(1-f)p$}
\\
(0,0) &\mbox{w. p. $1-p$.}
\end{cases}
\end{equation}

(Ep-2) If $(\Lf_n(\varepsilon), \Rf_n(\varepsilon))\in \{0\}\times \N$
\begin{equation}\label{TP21}
(\Lf_{n+1}(\varepsilon), \Rf_{n+1}(\varepsilon))=
\begin{cases}
(1,\Rf_n(\varepsilon)) \quad &\mbox{w. p. $fpr$}
\\
(0,\Rf_n(\varepsilon)+1) &\mbox{w. p. $(1-f)pr+p(1-r)$}
\\
(0,\Rf_n(\varepsilon)-1) &\mbox{w. p. $1-p$.}
\end{cases}
\end{equation}

(Ep-3) If $(\Lf_n(\varepsilon), \Rf_n(\varepsilon))\in \N\times \{0\}$
\begin{equation}\label{TP31}
(\Lf_{n+1}(\varepsilon), \Rf_{n+1}(\varepsilon))=
\begin{cases}
(\Lf_n(\varepsilon)+1,0) \quad &\mbox{w. p. $fpr+p(1-r)$}
\\
(\Lf_n(\varepsilon),1) &\mbox{w. p. $(1-f)pr$}
\\
(\Lf_n(\varepsilon)-1,0) &\mbox{w. p. $1-p$.}
\end{cases}
\end{equation}

(Ep-4) If $(\Lf_n(\varepsilon), \Rf_n(\varepsilon))\in \N\times \N$
\begin{equation}\label{TP41}
(\Lf_{n+1}(\varepsilon), \Rf_{n+1}(\varepsilon))=
\begin{cases}
(\Lf_n(\varepsilon)+1,\Rf_n(\varepsilon)) \quad &\mbox{w. p. $fpr+p(1-r)\varepsilon$}
\\
(\Lf_n(\varepsilon),\Rf_n(\varepsilon)+1) &\mbox{w. p. $(1-f)pr+p(1-r)(1-\varepsilon)$}
\\
(\Lf_n(\varepsilon)-1,\Rf_n(\varepsilon)) &\mbox{w. p. $1-p$.}
\end{cases}
\end{equation}
For $\varepsilon\in [0,1]$, we couple the processes $\{(\Lf_n(\varepsilon), \Rf_n(\varepsilon)): n \geq 1\}$ such that
\begin{align}\label{couple}
&\Lf_n (\varepsilon) \le \Lf_n (\varepsilon'),
\qquad\Rf_n (\varepsilon) \ge \Rf_n (\varepsilon')
\quad \mbox{ for }\varepsilon \le \varepsilon' \text{ and all } n \geq 1.
\end{align}

Taking $\Lf_n$, $\Rf_n$ and $N_n$ as in Subsection 2.1 and $\Lf_n (\cdot)$ and $\Rf_n (\cdot)$ as above, we have,  for
$\rho_n^f:=\frac{\Lf_n}{N_n}$, 
\begin{align}
&N_n(\varepsilon):=\Lf_n(\varepsilon)+\Rf_n(\varepsilon)=N_n \label{sumN}\\
& \Lf_{n+1}=\Lf_{n+1}\left(\rho_n^f\right),
\quad \Rf_{n+1}=\Rf_{n+1}\left(\rho_n^f\right),
\label{3:25}
\\
&\Lf_n(0) \le \Lf_{n} \le \Lf_n(1), \quad
\Rf_n(1) \le \Rf_{n} \le \Rf_n(0).
\label{3:27}
\end{align}
By the law of large numbers we have
\begin{align*}
\lim_{n\to\infty}\frac{L_n^f(\varepsilon)}{n}= \left[fpr+p(1-r)\varepsilon-1+p\right]_+,
 \text{ and } 
\lim_{n\to\infty}\frac{N_n}{n}= 2p-1, \text{ almost surely}, 
\end{align*}
and so, for $
\rho_n^f(\varepsilon):=\frac{\Lf_n(\varepsilon)}{N_n}
$, we have
\begin{align}
\lim_{n\to\infty}
\rho_n^f(\varepsilon)
&=\left[\frac{fpr+p(1-r)\varepsilon -1+p}{2p-1}
\right]_+
\nonumber
\\ \label{LLN}
&=\left[\frac{fpr-1+p}{2p-1}+ \frac{p(1-r)\varepsilon}{2p-1}\right]_+ .
\end{align}

We introduce the linear function defined by
$$
h(x)=\frac{fpr-1+p}{2p-1}+\frac{p(1-r)}{2p-1}x.
$$
Note that $\frac{p(1-r)}{2p-1}>0$.
By a simple calculation we see that if $f\le 1$
$$
h(0)\le \frac{pr-1+p}{2p-1}<0 \quad\text{ if $pr<1-p$}\; \quad\text{ and }\quad h(1)\le 1- \frac{2p}{2p-1} < 1.
$$
Then we may choose  $\delta>0$ such that
\begin{equation}\label{h<x}
\h(x) := h(x+\delta)< x, \quad [0,1].
\end{equation}
Put
\begin{align*}
\Lambda(\varepsilon, \delta)
&=\bigg\{ \omega :
\text{there exists } N=N(\omega)\in\N \mbox{ such that for all $n\ge N$, } 
\rho_n^f(\varepsilon) < \h(\varepsilon)
\bigg\}.
\end{align*}
From (\ref{LLN}), we have that 
\begin{align}\label{P=1}
P(\Lambda(\varepsilon, \delta))=1,
\quad \mbox{for all $\varepsilon, \delta\in (0,1]$.}
\end{align}
Also, 
taking $\varepsilon_c >0$ such that $\h(\varepsilon_c)=\eta$, i.e., 
$$
\varepsilon_c=\varepsilon_c(\delta,\eta):= \frac{\eta-\frac{fpr-1+p}{2p-1}}{\frac{p(1-r)}{2p-1}(1+\delta)}
=\frac{(2p-1)\eta-(fpr-1+p)}{p(1-r)(1+\delta)},
$$
we see that for $\varepsilon \le \varepsilon_c$ we have
$
\max\big\{\h(\varepsilon), \eta \big\}=\eta
$.

Now consider the recursion formula 
\begin{align}\label{R_F}
x_{n+1}= \h(x_n).
\end{align}
Since ($\ref{h<x}$), for $f<f_c<1$, 
\begin{equation}
\label{decreasing}
x_n \text{ is decreasing and }\lim_{n\to\infty}x_n= \frac{fpr-1+p}{pr-1+p}<0.
\end{equation}
We put
$$
\h(k,x):=\h(2^{-k} ([2^k x]+1)) \text{ for }k\in\N, 
$$
where $[a]$ the largest integer less than $a\in\R$.
From (\ref{decreasing}) we see that, for  sufficient large $k$, there exists $n_c\in \N$ such that 
\begin{align}\label{n_c}
\h^n(k,1):= \h(k,\h^{n-1}(k,1)) \le \eta \mbox{ for all } n \ge n_{c}.
\end{align}
Note that from (\ref{couple}) and (\ref{3:27}) we have that 
\begin{align}\label{inc}
\rho_n^f(\varepsilon) \ge \rho_n^f(\varepsilon')
 \mbox{ for } \varepsilon > \varepsilon' \text{, and } \rho_n^f \leq \rho_n^f(1),
\end{align}
thus, for any $\omega \in \bigcap_{m\in\N}\Lambda(m2^{-k}, \delta)$
there exists $N_1(\omega)$ such that,
for all $n\ge N_1(\omega)$,
$$
\rho_n^f[\omega] \le\rho_n^f(1)[\omega] \le \h(k,1), 
$$
and there exists $N_2(\omega)\ge N_1(\omega)$ such that
for all $n\ge N_2(\omega)$ 
$$
\rho_n^f[\omega] \le
\rho_n^f(\h(k,1))[\omega] \le \h(k,h(k,1))=\h^2(k,1).
$$
Repeating this procedure we have for any $\ell\in\N$ there exists $N_\ell(\omega)$ such that
for all $n\ge N_\ell(\omega)$ 
\begin{align}
\rho_n^f[\omega] \le \h^\ell (k,1).
\end{align}
From (\ref{n_c}), we now have
$$
\rho_n^f[\omega] \le \eta \quad \mbox{ for all $n\ge N_{n_c}(\omega)$}.
$$
Since $P(\bigcap_{m\in\N}\Lambda(m2^{-k}, \delta))=1$ from (\ref{P=1}), we have
\begin{align}
\lim_{n\to\infty}\rho_{n}^f[\omega] \le \eta,
\quad \mbox{a.s.}
\end{align}
Thus we obtain (\ref{key1}).

If
\begin{align}
\label{fcon}
f<  f_c -\frac{1-r}{r}\varepsilon,
\end{align}
$L_n^f(\varepsilon)$ is recurrent. Also, for $f < f_c\wedge 1$,  the condition (\ref{fcon}) holds for sufficiently small $\varepsilon$, hence
from (\ref{key1}) we see that $L_n^f$ hits the origin infinitely often.
This proves (i) of the Lemma \ref{lemma4}.

Let $f_c<1$. Observing that, for $\S_n^{f-}$ as in (\ref{r-bas}) and $\Lf_n(\cdot)$ as above,
$$
\S_n^{f-} \le \Lf_n(0),
$$
we see from (\ref{BAS1})-(\ref{BAS4}) that when $f>f_c$, for only finitely many $n$ we have $\S_n^{f-}=0$.
Thus, from  (\ref{3:27})  we have (ii) of (1).
Let $1 \le f_c < \frac{1}{r}$.
Since $1-p < p$ 
the random walk comparison as noted at the beginning of this section shows that $N_n \to \infty$ almost surely as $n \to \infty$. 
We have (ii) of (2).
\qed

\vspace{,5cm}

We give the proof of Theorem \ref{Th1}.
Part (1) is obtained by the random walk comparison. 
Part (3) is derived from (2) of Lemma 4.
The first statement of (2) is derived (ii) of (1) and (2) in Lemma \ref{lemma4}.

 Finally, considering the birth rate $rp$ of mutants, the limiting expected number of them with a fitness between $(a,b)$, with $f_c < a < b \leq 1$, is $rp(b-a)$. Thus we have, by an application of the strong law of large numbers
$$
\liminf_{n \to \infty} \frac{R_n^{b} - R_n^a}{N_n}  \geq \frac{p(b-a)}{2p-1} \text{ almost surely}.
$$
(Note this also follows from part (b) of the main Theorem of Guiol, Machado and Schinazi (2010).) 
This completes the proof of the second statement of part (2) of Theorem \ref{Th1}.

Finally, since the sites are each independently and uniformly distributed on $[0,1]$  Corollary \ref{GC} follows from Lemma \ref{lemma4}.

\section{Proof of Theorem \ref{Th2}}

We will prove Theorem \ref{Th2} with the help of two lemmas.

Let $A_k(t_1,n)$, $k, t_1,n\in\N$, be the event that a mutant born at time $t_1$ gets $k-1$ attachments until time $n$, and let  $q_k(t_1, n):=P(A_k(t_1,n))$. We have

\begin{lem}\label{Lemma6} Let $p=1$ i.e. no deaths.  For each $k,t_1 \in \N$
\begin{align}\label{429}
E\left[
\left\{
\frac{1}{n}\sum_{t_1=1}^n(\mathbf{1}_{A_k(t_1,n)}-q_k(t_1,n) )
\right\}^2
\right]\to 0 \text{ as }  n\to\infty.
\end{align}
\end{lem}
\noindent {\it Proof.}\quad 
The left hand side of (\ref{429}) is
\begin{align*}
&\frac{1}{n^2}\sum_{t_1=1}^n \sum_{s_1=1}^n \left[
P(A_k(s_1,n)\cap A_k(t_1,n)) - P(A_k(s_1,n))P(A_k(t_1,n))
\right]
\\
&=\frac{1}{n^2}\sum_{t_1=1}^n \sum_{s_1=1}^n P(A_k(s_1,n))\left[
P(A_k(t_1,n)\big{|}A_k(s_1,n)) - P(A_k(t_1,n))
\right].
\end{align*}
Thus it is enough to show the following for the proof of the lemma:
for any $x_1, y_1\in (0,1)$ with $x_1< y_1$
\begin{align}\label{430}
&P(A_k(y_1 n,n)\big{|}A_k(x_1 n,n))- P(A_k(y_1 n,n))\to 0, \quad n\to\infty.
\end{align}

Let $\{t_\ell\}_{\ell=1}^k$ be an increasing sequence of $\N$ with $t_k\le n$. 
We denote by $A_k[\{t_\ell \}_{\ell=1}^k;n]$  the event that a mutant comes at time $t_1$ which gets it's $(\ell-1)$th attachment at time $t_\ell$, $\ell= 2,3,\dots,k$, and no other attachment till time $n$.
Then
\begin{align}\label{sum_Ak}
A_k(t_1,n)=\sum_{\substack{t_2, t_3,\dots,t_k\in\N \\ t_1<t_2<\cdots < t_k<n}} A_k[\{t_\ell \}_{\ell=1}^k;n].
\end{align}
Let $\{s_\ell\}_{\ell=1}^k$  and $\{t_\ell\}_{\ell=1}^k$ be increasing sequences of $\N$ with $s_k,t_k\le n$. 

Suppose that $s_1=t_1$, then
\begin{align}\label{E1}
P(A_k[\{t_\ell \}_{\ell=1}^k;n]\big{|}A_k[\{s_\ell \}_{\ell=1}^k;n])= \bold{1}(s_\ell = t_\ell, \ell =2,3,\dots,k).
\end{align}
Also, for $s_1\not =t_1$,
if $\{s_\ell; \ell=2,\dots, k\}\cap \{t_\ell; \ell=2,\dots, k\}\not=\emptyset$, then
\begin{align}\label{E2}
P(A_k[\{t_\ell \}_{\ell=1}^k;n]\big{|}A_k[\{s_\ell \}_{\ell=1}^k;n])=0;
\end{align}
and if $\{s_\ell; \ell=1,2,\dots, k\}\cap \{t_\ell; \ell=1,2,\dots, k\}=\emptyset$, then
\begin{align*}
\nonumber
&P(A_k[\{t_\ell \}_{\ell=1}^k;n]\Big{|}A_k[\{s_\ell \}_{\ell=1}^k;n])
\\
\nonumber
&=P(A_k[\{t_\ell \}_{\ell=1}^k;n]\Big{|}\mbox{the mutant which came at time }t_1\\
\nonumber
& 
 \mbox{$\qquad\qquad\qquad\qquad\qquad$ does not get any attachment at times $\{s_\ell\}_{\ell=1}^k$ })
\\
&=P(A_k[\{t_\ell \}_{\ell=1}^k;n])\prod_{m : s_m >t_1}\left(1-\frac{\ell[s_m])(1-r)}{s_m} \right)^{-1},
\end{align*}
where $\ell[s_m]= \max\{ \ell : t_\ell < s_m\}$ is the population size at time $s_m$ of the fitness location occupied by
the mutant which came at time $t_1$.
Hence, we have,
\begin{align}\nonumber
&P(A_k[\{t_\ell \}_{\ell=1}^k;n]\big{|}A_k[\{s_\ell \}_{\ell=1}^k;n])- P(A_k[\{t_\ell \}_{\ell=1}^k;n])
\\ \nonumber
&= P(A_k[\{t_\ell \}_{\ell=1}^k;n]\big{|}A_k[\{s_\ell \}_{\ell=1}^k;n])
\left[
1-\prod_{m : s_m >t_1}\left(1-\frac{\ell[s_m])(1-r)}{s_m} \right)\right]
\\ \label{E3}
&\le \frac{k^2}{t_1} P(A_k[\{t_\ell \}_{\ell=1}^k;n]\big{|}A_k[\{s_\ell \}_{\ell=1}^k;n]).
\end{align}
Combining (\ref{E1}), (\ref{E2}) and (\ref{E3}) with (\ref{sum_Ak}),
we obtain (\ref{430}).
This completes the proof.
\qed

\vspace{.5cm}

Next we have
\begin{lem}\label{Lemma5} Let $p=1$. For each $k\in\N$
\begin{align}
&\lim_{n\to\infty}\frac{1}{n}\sum_{t_1=1}^n q_{k}(t_1,n)
=\frac{r}{1-r} B\left(\frac{2-r}{1-r}, k\right)=p_k.
\end{align}
\end{lem}

\noindent{\it Proof.}\quad 
Let $A_k(t_1,n)$  and   $q_k(t_1,n)$, $k, t_1,n\in\N$, be as above.
For $k=1$, we have
\begin{align*}
q_1(t_1, n)&= r \prod_{j=t_1+1}^n \left( 1-\frac{1-r}{j}\right),
\end{align*}
since the number of individuals at time $j-1$ is $j$ and the probability that the mutant who arrived at time $t_1$ gets an attachment at time $j$ is $\frac{1-r}{j}$. 

For $k=2$
\begin{align*}
q_2(t_1, n)&= r \sum_{t_2=t_1+1}^n \left\{\prod_{j=t_1+1}^{t_2-1} \left( 1-\frac{1-r}{j}\right)\right\}
\frac{1-r}{t_2}\left\{\prod_{j=t_2+1}^n \left( 1-\frac{2(1-r)}{j}\right)\right\},
\end{align*}
where $t_2$ is the time of the first attachment.
Similarly for each $k\in\N$
\begin{align*}
q_k(t_1, n)&= r \sum_{t_1 < t_2<\cdots <t_k\le n}
\prod_{\ell=1}^k \prod_{j=t_\ell +1}^{t_{\ell+1}}\left( 1-\frac{\ell(1-r)}{j}\right)
\prod_{\ell=1}^{k-1}\frac{\ell(1-r)}{t_{\ell+1}-\ell(1-r)},
\end{align*}
where we used the equation
$$
\frac{\ell(1-r)}{t_{\ell +1}}\frac{1}{1-\frac{\ell(1-r)}{t_{\ell +1}}}=\frac{\ell(1-r)}{t_{\ell+1}-\ell(1-r)}.
$$
By using Stirling's formula we see that
$$
\prod_{j=t_{\ell}+1}^{t_{\ell+1}} \left( 1-\frac{\ell(1-r)}{j}\right) \sim \left( \frac{t_{\ell}}{t_{\ell+1}}\right)^{\ell(1-r)},
\quad t_{\ell}, t_{\ell+1}\to\infty.
$$
Now letting  $n\to\infty$ and taking $t_\ell=n x_\ell$ we have
\begin{align*}
&\frac{1}{n}\sum_{t_1=1}^n q_{k}(t_1,n)
\sim r \int_{0<x_1<\cdots <x_{k}<1}dx_1\cdots dx_k \prod_{\ell=1}^k \left(\frac{x_\ell}{x_{\ell+1}}\right)^{\ell(1-r)}\prod_{\ell=1}^{k-1}\frac{\ell(1-r)}{x_{\ell+1}}
\\
&=r (1-r)^{k-1}(k-1)!\int_{0<x_1<\cdots <x_{k}<1}dx_1\cdots dx_k \;  x_1^{1-r} \prod_{\ell=2}^k x_{\ell}^{-r}
\\
&=r (1-r)^{k-1}\int_0^1 dx_1 x_1^{1-r} \prod_{\ell=2}^k \int_{x_1}^1 dx_\ell \; x_\ell^{-r}
\\
&= r \int_0^1 dx_1 x_1^{1-r} (1-x_1^{1-r})^{k-1}
\\
&=\frac{r}{1-r}\int_0^1 dy \; y^{\frac{1}{1-r}}(1-y)^{k-1}= \frac{r}{1-r} B\left(\frac{2-r}{1-r}, k\right).
\end{align*}
This compltes the proof. 
\qed

\vskip 3mm

We give the proof of Theorem \ref{Th2}.
When $p=1$ From Lemmas \ref{Lemma6} and \ref{Lemma5} we have
\begin{align*}
\frac{1}{n}\sum_{t_1=1}^n \mathbf{1}_{A_k(t_1,n)} 
 \to \frac{r}{1-r} B\left(\frac{2-r}{1-r}, k\right) \text{ as }  n\to\infty, 
\quad \text{in probability}.
\end{align*}
Noting that
$$
\lim_{n\to\infty}\frac{S_n}{n}=r, \quad \text{a.s.}
$$
we have
\begin{align}\label{con1}
&\lim_{n\to\infty}\frac{\sum_{f\in (0,1)}U_n^k(f)-U_n^k(f+)}{S_n}=\frac{1}{1-r} B\left(\frac{2-r}{1-r}, k\right)=p_k. \quad \text{in probability}.
\end{align}

Next we consider the case where $p\in (0,1)$.
We introduce another Markov process $\hat{X}_n$, $n\in \N\cup \{0\}$, which is a pure birth process, as follows:
\begin{enumerate}
\item At time $0$ there exists one individual at a site uniformly distributed on $(f_c,1)$.

\item with probability $p(1-rf_c)$ there is a new birth. There are two possibilities --
\begin{itemize}
\item with probability $\displaystyle{\hat{r}:=\frac{pr(1-f_c)}{p(1-rf_c)}}$ a mutant is born with a fitness uniformly distributed in  $[f_c,1]$,
\item with probability $\displaystyle{1-\hat{r}:=\frac{p(1-r)}{p(1-rf_c)}}$ a non-mutant individual is born. It has a fitness $f$ with a probability proportional to the number of individuals of fitness $f$,  and we increase the corresponding population of fitness $f$ individuals by $1$. 
\end{itemize}
\item With probability $1-p(1-rf_c)$  nothing happens, i.e. neither a birth nor a death occurs.
\end{enumerate}
For the Markov process $\hat{X}_n$, $n\in\N\cup \{0\}$,
we define $\hat{q}_k$, $\hat{S}_n$ and $\hat{U}_n$ in the same manner as $q_k$,  $S_n$ and $U_n$ for $X_n$, $n\in\N\cup \{0\}$. Then by the same argument as above we see that
\begin{align*}
\frac{1}{n}\sum_{t_1=1}^n \tilde{q}_{k}(t_1,n)&\sim p(1-rf_c) \frac{\hat{r}}{1-\hat{r}} B\left(\frac{2-\hat{r}}{1-\hat{r}}, k\right)
\end{align*}
and
$$
\lim_{n\to\infty}\frac{\hat{S}_n}{n}=pr(1-f_c).
$$
Hence
\begin{align*}
&\lim_{n\to\infty}\frac{\sum_{f\in (0,1)}\hat{U}_n^k(f)-\hat{U}_n^k(f+)}{\hat{S}_n}
=  \frac{1}{1-\hat{r}} B\left(\frac{2-\hat{r}}{1-\hat{r}}, k\right)
= p_k
\end{align*}
From Lemma \ref{lemma4}, we know that  deletions of individuals in $(f_c,1)$ occur finitely often and$\frac{R_n^f}{L_n^f+R_n^f}\to 1$ almost surely as $n\to\infty$. Thus we have 
\begin{align*}
&\lim_{n\to\infty}\frac{\sum_{f\in (0,1)}U_n^k(f)-U_n^k(f+)}{S_n}=\lim_{n\to\infty}\frac{\sum_{f\in (0,1)}\hat{U}_n^k(f)-\hat{U}_n^k(f+)}{\hat{S}_n} \quad \text{a.s.}
\end{align*}
and so (\ref{con1}) for $p\in (0,1]$.
Noting that the sites are uniformly distributed on $[0,1]$ independently, and preferential attachment does not depend on the position of sites, we obtain Theorem \ref{Th2} from (\ref{con1}).
\qed

\vspace{.6cm}
\section{Number of individuals of a fixed fitness}
Fix $f\in [0,1]$ and let $N_n^f$ denote the number of individuals with fitness $f$ at time $n$. When $rp > 1-p$, i.e. $f_c < 1$,  from Lemma \ref{lemma4} we know that, $P(L_n^f=0 \text{ infinitely often}) =1$ for $f \in (f_c,1)$. Thus, if a mutant with fitness $f \in (f_c,1)$ is born at some large time $\ell$, then the chances of the mutant dying is small, and so a natural question is `for some $n > \ell$, how many individuals did this mutant attract by  time $n$', i.e., what is the value of $N_n^f$?
\begin{prop}
Fix $f\in (f_c,1)$, we have, for $\ell < n$,  as $ \ell, n\to\infty$
\begin{align*}
&E[N_n^f | \text{a mutant with fitness $f$ is born at time $\ell$}]\\
&\sim \frac{\Gamma((2p-1)\ell+1)\Gamma((2p-1)n+1+p(1-r))}{\Gamma((2p-1)\ell+1+p(1-r))\Gamma((2p-1)n+1)} 
\\
&\sim\left( \frac{n}{\ell}\right)^{p(1-r)}.
\end{align*}
\end{prop}

\noindent {\it Proof}. Since we are interested in the region $f > f_c$ and also, for the calculation of the expectation, we just need to factor out the death rate $(1-p)$, so we modify the Markov process  
$\hat{X}_n$ introduced in the proof of Lemma \ref{Lemma5}, by removing the times when `nothing happens' , i.e. the process does not move. This is done as follows: let $\hat{N}_n$ be the number of individuals of  the process $\hat{X}_n$ at time $n$, we define a new Markov process $\check{X}_n$, for $n \geq 0$, by
$$
\hat{X}_n = \check{X}_{\hat{N}_n-1}.
$$
Since $\hat{N}_0 = 1$, we see that $\check{N}_\ell = \ell+1$, where $\check{N}_\ell$  is the number of individuals of the process  $\check{X}$ at time $\ell$.

Letting $\check{N}_m^f$ denote the number of individuals of the $\check{X}$ process of fitness $f$ at time $m$, we have
\begin{align*}
&E[\check{N}_m^f|\check{N}_{m-1}^f]
\\
&=\{1-p(1-r)\}\check{N}_{m-1}^f +p(1-r)\left\{(\check{N}_{m-1}^f+1)\frac{\check{N}_{m-1}^f}{m}+\check{N}_{m-1}^f \left(1-\frac{\check{N}_{m-1}^f}{m} \right)\right\}
\\
&=\left(1+\frac{p(1-r)}{m}\right)\check{N}_{m-1}^f.
\end{align*}
If $\check{N}_0^f=\check{N}_0=1$ then we have
\begin{align}
&E[\check{N}_m^f | \check{N}_0^f=1]=\prod_{k=1}^{m}(\frac{k+p(1-r)}{k})=\frac{\Gamma(m+1+p(1-r))}{\Gamma(1+p(1-r))\Gamma(m+1)},
\end{align}
while, if $\check{N}_\ell^f=1$ then we have
\begin{align}
&E[\check{N}_m^f | \check{N}_\ell^f=1]= \prod_{k=\ell+1}^{m}(\frac{k+p(1-r)}{k})
=\frac{\Gamma(\ell+1)\Gamma(m+1+p(1-r))}{\Gamma(\ell+1+p(1-r))\Gamma(m+1)}.
\end{align}
Since $\frac{\hat{N}_n}{n}\to pr(1-f_c)+p(1-r)=2p-1$, 
if $\hat{N}_0^f=1$ then we have
\begin{align}
&E[\hat{N}_n^f | \hat{N}_0^f=1]\sim \prod_{k=1}^{(2p-1)n}(\frac{k+p(1-r)}{k})=\frac{\Gamma((2p-1)n+1+p(1-r))}{\Gamma(1+p(1-r))\Gamma((2p-1)n+1)}.
\end{align}
Also,  $\frac{\hat{N}_\ell}{\ell}\to 2p-1$, so
for $\hat{N}_\ell^f=1$, we have
\begin{align*}
E[\hat{N}_n^f | \hat{N}_\ell^f=1]&= \prod_{k=(2p-1)\ell+1}^{(2p-1)n}(\frac{k+p(1-r)}{k})
\\
&=\frac{\Gamma((2p-1)\ell+1)\Gamma((2p-1)n+1+p(1-r))}{\Gamma((2p-1)\ell+1+p(1-r))\Gamma((2p-1)n+1)}.
\end{align*}
From Lemma \ref{lemma4} we have $E[N_n^f |N_\ell^f=1]\sim E[\hat{N}_n^f|\hat{N}_\ell^f=1]$,
and that completes the proof of the proposition.
\qed
\section{Heuristics for the case $f_c > 1$}
We now present some mean field heuristics about the location of the leftmost site
 $x_t$ at time $t$ in the case when $pr < 1-p < p$, i.e. $f_c > 1$. These heuristics should be seen in connection with (\ref{Thm11}) of Theorem \ref{Th1}.
 
 Let $y_t=1-x_t$. The number of individuals to enter the interval $(x_t,1]$ is approximately
$$
pry_t dt + p(1-r)dt, 
$$
where the first term counts the births which are mutants and the second term counts the births which are not mutants. While the number of individuals deleted in the interval $(x_t,1]$ is approximately
$$
-\frac{dy_t}{y_t} \{ p-(1-p) \}t,
$$
this being  the absolute value of the deletions since $\frac{dy_t}{dt}<0$.
Thus we consider the following differential equation:
\begin{align*}
pry_t dt + p(1-r)dt +\frac{dy_t}{y_t}\{ 2p-1\}t=(2p-1)dt,
\end{align*}
from which we have
\begin{align*}
\frac{dt}{t} &= \frac{-(2p-1)}{pry_t+p(1-r)-(2p-1)}\frac{dy_t}{y_t}
\\
&= \frac{-(2p-1)}{pry_t+pr(f_c-1)}\frac{dy_t}{y_t}
\\
&= -\frac{2p-1}{pr}\left\{ \frac{1}{y_t+ (f_c-1)} \right\}\frac{dy_t}{y_t}
\\
&=-\frac{2p-1}{pr(f_c-1)}\left\{ \frac{f_c-1}{y_t+ (f_c-1)} \right\}\frac{dy_t}{y_t}
\\
&=-\frac{2p-1}{pr(f_c-1)}\left\{ 1 - \frac{y_t}{y_t+ (f_c-1)} \right\}\frac{dy_t}{y_t}
\\
&=-\frac{2p-1}{pr(f_c-1)}\left\{ \frac{1}{y_t} - \frac{1}{y_t+ (f_c-1)} \right\}dy_t.
\end{align*}
Hence, for an appropriate constant $c$, we have
$$
c+ \log t = \frac{2p-1}{pr(f_c-1)}\left[\log (y_t+(f_c-1)) -\log y_t\right]=  \frac{2p-1}{pr(f_c-1)}\log(1+ \frac{1-r}{y_t}),
$$
and so
\begin{align*}
&1+ \frac{f_c-1}{y_t}=\exp\left\{\{c+\log t)\frac{pr(f_c-1)}{2p-1}\right\}=Ct^{\gamma}, 
\end{align*}
where $\displaystyle{ \gamma = \frac{pr(f_c-1)}{2p-1} = \frac{1-p-pr}{2p-1}}$
and $C=e^{c\gamma}$. Thus
$$
y_t = \frac{f_c-1}{Ct^\gamma -1} \sim C' t^{-\gamma}, \quad t\to\infty.
$$
Moreover, the number of sites is approximately
$$
rpty_t\sim C'r(1-r)p t^{1-\gamma}.
$$
\begin{remark}
For $f_c >1$ we have $\gamma=\gamma(p,r) > 0$, and $\gamma(p,r)$ is a decreasing function of $p$. Also
\begin{itemize}
\item [(i)] when  $p=1-p$, i.e., $p=\frac{1}{2}$, then $\gamma = \infty$; this corresponds to the case when the process dies out repeatedly,
\item[(ii)] when $pr=1-p$, i.e., $f_c=1$, then $\gamma = 0$; this corresponds to the case when the number of sites surviving is of order $o(t)$.
\item[(iii)] when $p=\frac{2}{3+r}\in \left(\frac{1}{2},\frac{1}{1+r} \right) $, then $\gamma =1$; this corresponds to the case when there are only a bounded number of sites surviving.
\end{itemize}
\end{remark}
From the above, we see that there are three critical values
$$
p_c^{(0)}:=\frac{1}{2} < p_c^{(1)}:= \frac{2}{3+r} < p_c^{(2)}:=\frac{1}{r+1}<1
$$
and four phases:
\begin{enumerate}
\item For $p\in (p_c^{(2)}, 1)$,  $\gamma \in (-r,0)$  and individuals exist in the interval $(f_c,1]$.

\item For $p\in (p_c^{(1)}, p_c^{(2)})$, $\gamma \in (0,1)$ and the number of sites are increasing with the order $t^{1-\gamma}$
and the average number of particles is of order $t^{\gamma}$.

\item For $p\in (p_c^{(0)}, p_c^{(1)}]$,  $\gamma \in (1,\infty)$, that is, $1-\gamma$ is negative, and the number of sites is finite, with the average number of individuals being of order $t$.

\item For $p\in (0,p_c^{(0)}]$ the process dies out infinitely often.

\end{enumerate}

\section{Simulation}
We conclude the paper with some simulations. The R code is given in the appendix. Here we have taken $p=3/4$, $r=1/2$, so that $f_c = 2/3$. The simulation has been conducted with $n = 100,000$.

\begin{figure}[htp] \centering{
\includegraphics[scale=0.39]{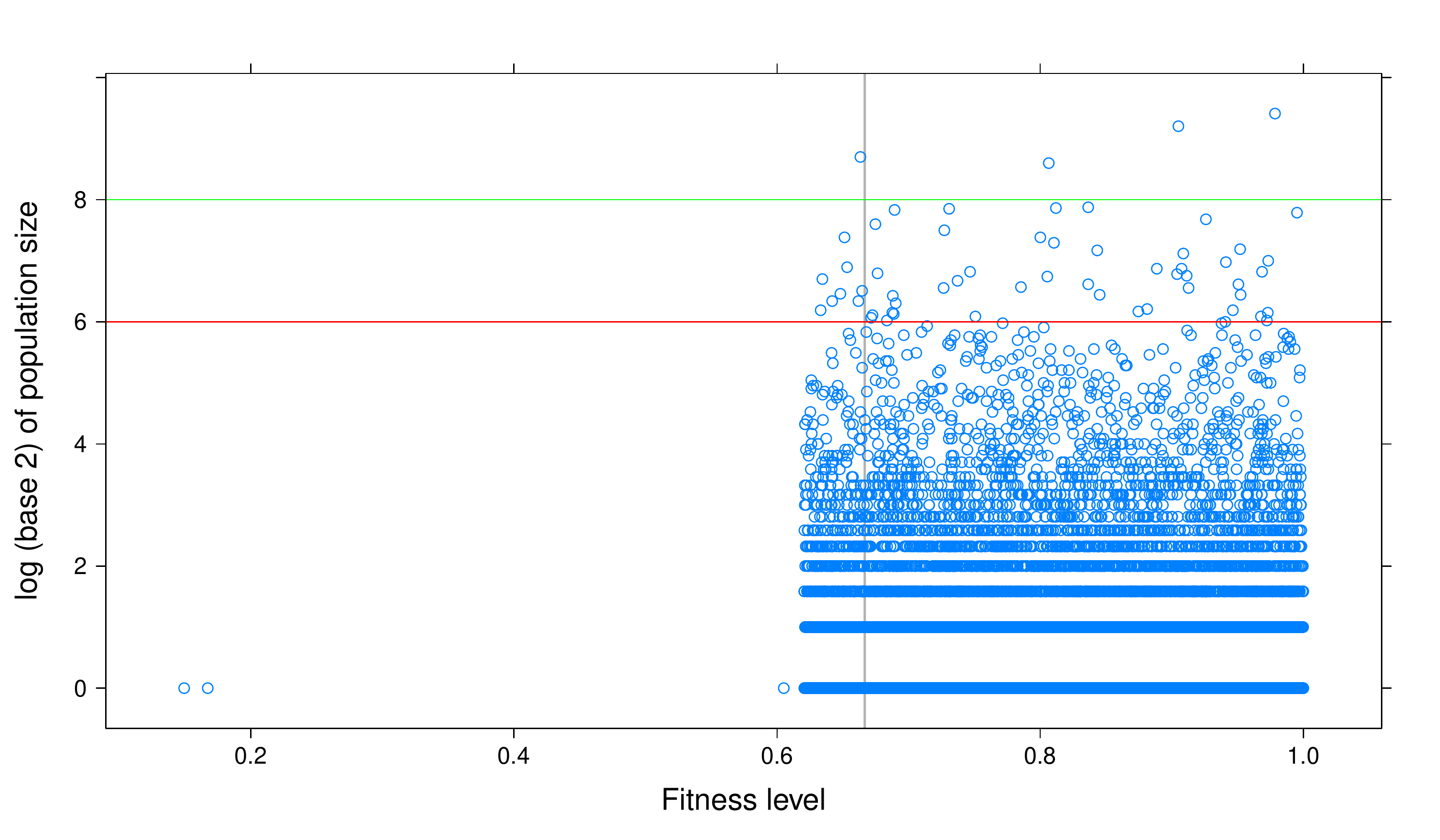}
}
\caption{Population (in $\log_2$ scale) at various fitness levels.}
\label{popsize-fitness}
\end{figure}

Figure \ref{popsize-fitness} presents the size of the population in $\log_2$ scale at each surviving site. The plot above the red line indicates the sites where the population size is $2^6$ or more, while the plot above the green line indicates the sites where the population size is $2^8$ or more.

\begin{figure}[htp] \centering{
\includegraphics[scale=0.35]{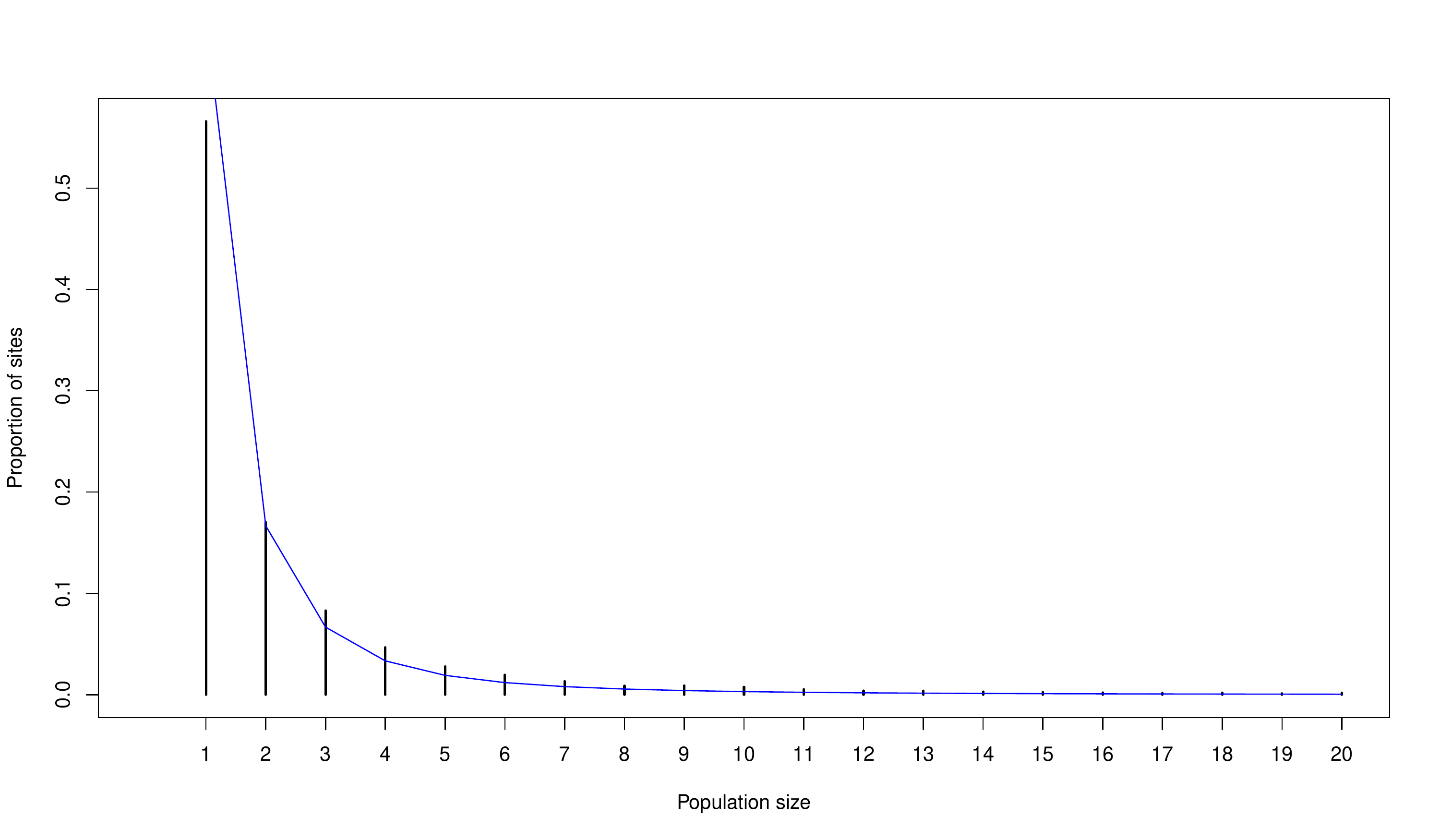}
}
\caption{Theoretical and observed proportion of sites with respect to population size.}
\label{kdist}
\end{figure}

In Figure \ref{kdist} the $x$-axis gives the population size, while the y-axis presents the proportion of sites with the given population size. The blue line is the theoretical value as obtained from Theorem \ref{Th2} and the vertical bars are the observed values.

\section{Acknowledgements}
The authors are grateful to Professor Deepayan Sarkar who wrote the R code and performed the simulation.
Rahul Roy acknowledges the grant MTR/2017/000141 from DST which supported this research and also the hospitality of Keio University where much of the work was done.
Hideki Tanemura's research is  supported in part by Grant-in-Aid for Scientific Research (S), No.16H06338; Grant-in-Aid for Scientific Research (B), No.19H01793 from Japan Society for the Promotion of Science.

\vspace{.6cm}

\noindent{\bf Appendix: The R code for the simulation}

\begin{verbatim}

library(lattice)
library(latticeExtra)


createState <- function(MAX_POP = 10000L, p = 3/4, r = 1/2)
{
    n <- integer(MAX_POP)   # size of each sub-population
    f <- numeric(MAX_POP)   # fitness of each sub-population
    tob <- integer(MAX_POP) # time at which this population first appeared
    n[1] <- 1L
    f[1] <- 0
    npop <- 1L
    ndead <- 0L
    t <- 0L
    environment()v
}
## Make sure to keep normalized by ordering f from low to high
updateState <- function(S)
{
    p <- S$p
    r <- S$r
    S$t[] <- S$t + 1L # increment process lifetime counter
    u <- runif(1) # to decide which branch
    f <- runif(1) # new fitness value if needed
    if (u < 1-p) # kill particle with lowest fitness
    {
        if (S$n[1] > 0L) S$n[1] <- S$n[1] - 1L
        if (S$npop > 0 && S$n[1] == 0L) { # a population has just died out
            S$ndead[] <- S$ndead + 1L
            S$n[1:S$npop] <- S$n[2:(S$npop+1)]
            S$f[1:S$npop] <- S$f[2:(S$npop+1)]
            S$tob[1:S$npop] <- S$tob[2:(S$npop+1)]
            S$npop[] <- S$npop - 1L
        }
    }
    else if (u < 1 - p + p * r || S$npop == 0) # create new sub-population\\
    {
        S$npop[] <- S$npop + 1L
        if (S$npop == S$MAX_POP)
            stop("exceeded maximum sub-populations allowed: ", S$MAX_POP)
        S$f[S$npop] <- f
        S$n[S$npop] <- 1L
        S$tob[S$npop] <- S$t
        i <- 1:S$npop
        ord <- order(S$f[i])
        S$n[i] <- (S$n[i])[ord]
        S$f[i] <- (S$f[i])[ord]
        S$tob[i] <- (S$tob[i])[ord]
    }
    else { # increment size of one population by 1
        i <- sample(S$npop, 1, prob = S$n[1:S$npop])
        S$n[i] <- S$n[i] + 1L
    }
}


S <- createState(MAX_POP = 20000, p = 3/4, r = 1/2)
(f_c <- with(S, (1-p) / (p*r)))


for (i in 1:100000) updateState(S)

Sdf <- subset(as.data.frame(as.list(S)), n > 0, select = c(tob, f, n))
names(Sdf) <- c("time of birth", "fitness", "population size")
xyplot(log2(`population size`) ~ fitness, data = Sdf, cex = 0.7, 
       ylab = "log (base 2) of population size", xlab = "Fitness level",
       abline = list(v = f_c, col = "grey70", lwd = 2)) + layer(panel.abline(h = c(6, 8), col = c("red", "green")))

pk <- function(k, r) 1 / (1-r) * beta((2-r) / (1-r), k)
plot(prop.table(table( Sdf[["population size"]] )), 
     xlim = c(0, 20), xlab = "Population size", ylab = "Proportion of sites")
lines(1:20, pk(1:20, r = S$r), col = "blue")

\end{verbatim}

\vspace{.5cm}
\noindent Rahul Roy\\
Theoretical Statistics and Mathematics Unit\\
Indian Statistical Institute\\
7 SJS Sansanwal Marg\\
New Delhi 110016, India.\\

\vspace{.3cm}
\noindent Hideki Tanemura\\
Department of Mathematics \\
Keio University\\
 Hiyoshi, Kohoku-ku\\
 Yokohama 2238522, Japan 


\end{document}